\documentclass[pdflatex,sn-mathphys]{sn-jnl}
\usepackage{enumitem}
\DeclareMathAlphabet{\mathpzc}{OT1}{pzc}{m}{it}

\numberwithin{equation}{section}

\newcommand\isoto{\stackrel{\sim}{\smash{\longrightarrow}\rule{0pt}{0.4ex}}}
\newcommand\isomto{\stackrel{\textstyle\sim}{\smash{\longrightarrow}\rule{0pt}{0.4ex}}}

\theoremstyle{thmstyleone}
\newtheorem{theorem}{Theorem}[section]
\newtheorem{proposition}[theorem]{Proposition}%

\theoremstyle{thmstyletwo}%
\newtheorem{example}{Example}[section]%
\newtheorem{remark}{Remark}%

\theoremstyle{thmstylethree}%
\newtheorem{definition}{Definition}[section]%

\begin{document}

\title[An Infinitary Model of the Diagrammatic Calculus]{An Infinitary Model of Diagrammatic Calculus in Unbiased Monoidal Categories}
\author{\fnm{Vihaan} \sur{Dheer}}

\abstract{Properties of morphisms represented by so-called `string diagrams' of monoidal categories (and their braided and symmetric derivatives), mainly their resistance in value to isotopic deformation, have made the usage of graphical calculi commonplace in category theory ever since the correspondence between diagrams and tensor categories was rigorously established by Joyal and Street in 1991. However, we find it important to make certain additions to the existing theory of monoidal categories and their diagrams, with the goal of extending to so-called `infinitary monoidal categories'. Most crucially, we employ a structure inherently resistant to isotopic deformation, thus replacing topological details with categorical ones. In the process, we coherently introduce infinitary tensor product and transfinite composition into the diagrammatic formalism.}


\maketitle

\tableofcontents

\section{Introduction}\label{intro}
\subsection{Motivations}
The first successful design and application of the concept of string diagrams similar to those employed by monoidal category theorists was likely due to Roger Penrose in 1971, coining \textit{Penrose notation} for tensors, in a series of works beginning with Ref. \citenum{penrose-1}, featuring various others including Ref. \citenum{penrose-2}. In 1991, the celebrated work of Joyal and Street ``The Geometry of Tensor Calculus, I'' \cite{main-source-js} formally introduced this concept into the scope of monoidal category theory. Ever since MacLane's definition of the monoidal (tensor) category \cite{categories-smc}, much effort had been placed in the study of the coherence theorem \cite{nat-assoc}, as well as the reformulation of the axioms \cite{kelly-conditions}. Upon Joyal and Street's application of string diagrams to the theory, addressing not only standard monoids but braided and symmetric categories as well, attention was diverted towards studying and re-purposing these diagrams to make use of in related categories. Twenty years from then, it had become commonplace in monoidal category theory to such a degree as to prompt expository works like Ref. \citenum{survey-graph-types}.

More recently, as monoidal categories have begun to gain applicability to external fields such as quantum information theory \cite{quantum-ct} and topological field theory \cite{tqft-reference}, and their diagrams as a result have become more widely applied, the graphical theory may require modifications to stay in use. In reference specifically to the proof for Theorem 1.2 of Ref. \citenum{main-source-js}, we aim to make a few major updates, from which a new diagrammatic theory supportive of more general kinds of diagram will arise. 

The mentioned theorem states that isotopic deformation of \textit{progressive plane} (string) diagrams does not change its `value', i.e. the morphism which it depicts. Clearly this makes checking validity of equations in these categories much easier, but it is only half of the full picture: one may not only be interested in isotopy-equivalence implying equality, but also whether equality of morphisms implies that their diagrams can be isotopically deformed into one another. Subtleties in this question, related to the various representations of a single morphism, are valuable to consider to the goal of furthering the theory.

Furthermore, Joyal and Street's notion of progressive plane diagrams do not allow for so-called `infinitary diagrams', a concept which the author has yet to see mentioned in the relevant literature, but feels it has potential for later use. The reason is likely as follows: string diagrams have historically been primarily for the purpose of \textit{human} computation, thus the difficulties in visualizing infinite diagrams (in either dimension) may have swayed most away from the concept. However, the \textit{monoidal diagrams} we introduce in this work refer not just to the visual concept, but are a new type of mathematical structure, which have use in their own right, not just as visual aids.

Our method for making these adjustments in some ways \textit{categorifies} the original proof; in other words, at least initially, diagrams are not topological entities, and in fact are by design \textit{inherently} resistant to isotopic deformation, a concept which will be made formal in later sections. Instead, they are mathematical structures which shall become members of a `category of diagrams' induced by some monoidal category; in the finitary case, the two formulations are of course equivalent.

To support infinite arity of diagrams, we will rely heavily on Tom Leinster's \cite{unbiased-mc} \textit{unbiased} definition of monoidal categories, and alter it to define the notion of $\alpha$-ary monoidal categories, allowing us to work with infinite monoidal products. Finally, we shall employ the standard notion of transfinite composition in categories, effectively permitting categories to be of infinite length in both the horizontal and vertical directions.

\subsection{Outline}
This paper will be structured as follows: in Sec. \ref{cat-prelim}, we will review transfinite composition, introduce our generalized definition of monoidal categories, and discuss a novel construction called the `colimit expansion'. Sec. \ref{monoidal-diagrams} firstly defines \textit{monoidal diagrams}, and then focuses on extracting from such diagrams the morphisms in the underlying category that they encode. Finally, we conclude with a discussion of applications and future work in Sec. \ref{conclusion}.

\section{Infinitary Monoidal Categories}\label{cat-prelim}
This section will develop the categorical formalism necessary to extend standard monoidal categories to their `infinitary counterparts', in which infinite products and compositions exist. This is important both to study the monoidal diagrams that will be discussed in the next section, but also because these categories are interesting in their own right. After reviewing transfinite composition, we define these infinitary categories by converting the axioms for a standard monoidal category into their unbiased form, and then extending via ordinals. Finally, we conclude with a discussion of the so-called ``colimit-expansion", which provides a systematic method to extend standard monoidal categories to infinitary ones under certain conditions.
\subsection{Transfinite Composition}
Extending monoidal categories ``horizontally" will involve considering compositions of infinite sequences. Such `transfinite compositions' appear all over mathematics, perhaps most commonly in the study of model categories (see e.g. Ref. \citenum{model-categories}). We briefly define transfinite composition and discuss some basic results relevant here.
\begin{definition}[transfinite composition]\label{def:transfinite-composition}
In a category $\mathscr{C}$, for some ordinal $\alpha$ (taken as the corresponding thin category), an $\alpha$-\textit{transfinite sequence} is a functor
\begin{equation}
    X_\bullet\colon\alpha\to\mathscr{C}
\end{equation}
such that $X_\beta=\varinjlim_{\gamma<\beta}X_\gamma$ for each non-zero limit ordinal $\beta<\alpha$. The \textit{transfinite composition} of such a sequence is the unique morphism
\begin{equation}
    X_0\longrightarrow \varinjlim_{\gamma<\alpha}X_\gamma
\end{equation}
induced by the colimit. The existence of this colimit is guaranteed by transfinite induction.
\end{definition}
The most intuitive example is given by successive inclusions in a concrete category $(\mathscr{C},F\colon\mathscr{C}\to\mathbf{Set})$, for $F$ faithful. A sequence
\begin{equation}
\begin{tikzcd}
X_0 \arrow[r, hook] & X_1 \arrow[r, hook] & X_2 \arrow[r, hook] & \cdots
\end{tikzcd}
\end{equation}
of inclusions under $F$ has a transfinite composition that is the obvious inclusion $X_0\subseteq\varinjlim X_\bullet$ of underlying sets.

A basic but useful result is that two transfinite compositions are compatible (i.e. rearrangeable) with each other:
\begin{proposition}
Let $X_\bullet,X'_\bullet$ be $\alpha,\alpha'$-transfinite sequences respectively such that $X'_0=\varinjlim X$. Let $\beta,\beta'$ be ordinals such that $\alpha+\alpha'=\beta+\beta'$; then we may form $\beta,\beta'$-transfinite sequences $Y,Y'$ respectively which are segments/concatenations of $X,X'$. Then $Y_0=X_0$, $Y'_0=\varinjlim Y$, and $\varinjlim Y'=\varinjlim X'$.
\end{proposition}
Under certain conditions, a class of morphisms in a category may be closed under transfinite composition; in other words the the composition of any transfinite sequence with contained morphisms exists and is as well contained. We call such a class \textit{saturated}. As standard compositions are $2$-transfinite, a saturated class is necessarily a subcategory. One example is that of monomorphisms.
\begin{proposition}
The wide subcategory of monomorphisms $\mathrm{Mono}(\mathscr{C})$ is saturated.
\begin{proof}
We must show that for any transfinite sequence of monomorphisms
\begin{equation}
    X_\bullet\colon\alpha\hookrightarrow\mathscr{C},
\end{equation}
the canonical map $X_0\longrightarrow\varinjlim_\gamma X_\gamma$ is monic. Suppose it weren't; by transfinite induction, no leg of the limiting cocone is monic. But this contradicts the uniqueness given by the universal property of $\varinjlim_\gamma X_\gamma$.
\end{proof}
\end{proposition}
In fact, this is a special case of a more general result stating that the morphism class with the `left lifting property' with respect to any class of morphisms is saturated \cite{higher-cat-homotopical-alg}. This of course has various corollaries; for example, the classes of cofibrations and trivial cofibrations in a model category are saturated.

In the following discussion, when making categories compositionally infinitary, we will use sequential cocompleteness to ensure a category has all transfinite compositions. These saturated subcategories are especially well-behaved because transfinite compositions in the subcategory coincide with those in the category itself.

\subsection{Unbiased Monoidal Categories}\label{sec-def-unb}
The contrast between biased and unbiased definitions is ubiquitous in mathematics. In short, we are often given an operation intended to be of arity $n$, that is presented solely in the nullary and binary cases. In this case the definition is called \textit{biased}, contrary to an \textit{unbiased} definition which constructs an $n$-ary operation at the outset. 

In no way is the term meant to suggest some superiority of unbiased definition: many cases, for example the monoidal categories at hand, require much shorter and more efficient definitions, as well as auxiliary coherence conditions. However, generalization with ordinal numbers is much easier with an unbiased definition, which is the motivation for using one here.

We first formally define partitions of ordinal numbers, which greatly shrinks the coherence conditions necessary for unbiased monoidal categories, which shall be introduced thereafter. Finally, we will redevelop some of the relevant theory of standard monoidal categories for their unbiased variants.
\begin{definition}[ordinal partitions]
A \textit{$\gamma$-partition} of an ordinal $\alpha$ is a $\gamma$-sequence which sums to alpha. Explicitly, it is a function $p_\bullet\colon\gamma\to\alpha$ such that
\begin{equation}
    \sum_{i<\gamma}{p_i}=\alpha.
\end{equation}
The set of all $\gamma$-partitions of $\alpha$ is written $\alpha^\gamma_\Sigma$, while the set of all partitions in general of $\alpha$ is $\alpha^<_\Sigma$. \textit{A $(\delta_\bullet,p_\bullet)$-double partition of $\alpha$} is an indexed family $p_{\bullet, \bullet}=\{p_{i,\bullet}\in (p_i)_{\Sigma}^{\delta_i}\,\colon i<\gamma\}$, for $\gamma$-partitions $\delta_\bullet, p_\bullet$ of ordinals $\delta, \alpha$ respectively. For each such double partition there is an \textit{induced $\delta$-partition of $\alpha$} $(p_{\bullet,\bullet})^\cup_\bullet\colon\delta\to\alpha$, or simply $(p^\cup)_\bullet$, defined in the obvious fashion. When the domain and range of a partition are clear, then an ordinal $\beta$ acts as a \textit{constant} partition $\beta\colon\gamma\to\delta$ such that $i\mapsto\beta$.
\end{definition}

\begin{definition}[$\alpha$-ary unbiased monoidal category]\label{unb-mon}
An \textit{$\alpha$-ary unbiased monoidal category} for some ordinal $\alpha\nless\omega$ (where $\alpha=\mathbf{Ord}$ is permissable) is a tuple $(\mathscr{C}, \mathcal{T}_\bullet, \iota, \eta)$, with an ordinary category $\mathscr{C}$, a formal family of functors
\begin{equation}
    \mathcal{T}_\bullet\colon\alpha\to\bigcup_{\beta<\alpha}{\mathrm{Hom}{\left(\mathscr{C}^\beta,\mathscr{C}\right)}}\quad\text{s.t.}\quad \mathcal{T}_\beta\colon\mathscr{C}^\beta\to\mathscr{C},
\end{equation}
a natural isomorphism $\iota:\mathrm{id}_{\mathscr{C}}\isomto\mathcal{T}_1$, and a family of natural isomorphisms
\begin{equation}
    \eta=\big\{\eta_{p_\bullet}\colon \mathcal{T}_\gamma\circ\prod_{i<\gamma}{\mathcal{T}_{p_i}}\isomto\mathcal{T}_\beta \,\big\vert\, p_\bullet\in\beta^\gamma_\Sigma\,\,\,\forall\gamma,\beta<\alpha \big\},
\end{equation}
such that for any $\gamma$-partition $p_\bullet$ of $\beta<\alpha$, and any $(\delta_\bullet, p_\bullet)$-double partition, the following diagram commutes,
\begin{equation}\label{eq:unb-ax-1}
\begin{tikzcd}
{\mathcal{T}_\gamma\circ\prod_{i<\gamma}{\left(\mathcal{T}_{\delta_i}\circ\prod_{j<\delta_i}{\mathcal{T}_{p_{i,j}}}\right)}} \arrow[rrr, "{\mathcal{T}_\gamma\left(\prod_{i<\gamma}{\eta({p_{i,\bullet}})}\right)}"] \arrow[dd] &  &  & \mathcal{T}_\gamma\circ\prod_{i<\gamma}{\mathcal{T}_{p_i}} \arrow[dd, "\eta(p_\bullet)"] \\
&  &  &                                                              \\
\mathcal{T}_\delta\circ\prod_{k<\delta}{\mathcal{T}_{(p^\cup)_k}} \arrow[rrr, "\eta\big((p^\cup)_\bullet\big)"]                        &  &  & \mathcal{T}_\beta                                            
\end{tikzcd}
\end{equation}
as does the following:
\begin{equation}\label{eq:unb-ax-2}
\begin{tikzcd}
\mathcal{T}_\beta \arrow[rrddd, "\mathrm{id}_{\mathcal{T}_\beta}"] \arrow[ddd, "\iota\mathcal{T}_\beta"] \arrow[rr, "\mathcal{T}_\beta\left(\prod_{i<\gamma}{\iota}\right)"] &  & \mathcal{T}_\beta\circ\prod_{i<\gamma}{\mathcal{T}_1} \arrow[ddd, "\eta\left(1_\bullet\right)"] \\
&                                                                                                 \\
&                                                                                                 \\
\mathcal{T}_1\circ\mathcal{T}_\beta \arrow[rr, "\eta(\beta_\bullet)"]             &  & \mathcal{T}_\beta                                                            
\end{tikzcd}
\end{equation}
\end{definition}
\begin{definition}[infinitary, finitary, strict]
An $\alpha$-ary unbiased monoidal category is called \textit{finitary} if $\alpha=\omega$, and \textit{infinitary} otherwise. If $\alpha=\mathbf{Ord}$, it is \textit{fully infinitary}. An $\alpha$-ary unbiased monoidal category is \textit{strict} if each component $\iota$ of each element of $\eta$ in Definition \ref{unb-mon} is an identity.
\end{definition}
From any unbiased monoidal category we can build a biased monoidal category by extracting the associator/unitors from the family $\eta$ (see the proof of Thm. \ref{thm:strictification}). Likewise, given a biased monoidal category, we can build a finitary unbiased monoidal category; the two notions are thus equivalent in the finitary case.
\begin{remark}
In the definition, we have made implicit reference to the categories $\mathscr{C}^\beta$ for an ordinal $\beta$. To make this notion complete, we wrap these objects in a \textit{power category of $\mathscr{C}$} which includes all ordinal powers of $\mathscr{C}$ and the functors between them as morphisms. This category becomes the prototypical example of an unbiased strict (skeletal) monoidal category when endowed with the Cartesian product, producing objects for general ordinal sequences. The strictness derives from the fact that we do \textit{not} distinguish between $\mathscr{C}^\beta\times\mathscr{C}^\gamma$ and $\mathscr{C}^{\beta+\gamma}$.
\end{remark}
The most obvious example is given by generalizing (co)cartesian categories. Any category $\mathscr{C}$ with all small (co)products will be fully infinitary if each $\mathcal{T}_\gamma$ is taken to be the generalized (co)product. A less trivial example is provided by the category of endofunctors (that of which a monad is an object).
\begin{example}
Suppose that, for a small category $\mathscr{C}$, the full subcategory $\{\mathscr{C}\}\subseteq\mathbf{Cat}$ is cocomplete (or, in principle, sequentially cocomplete). Then transfinite compositions of endofunctors in the sense of Def. \ref{def:transfinite-composition} form a fully infinitary product on the strict monoidal category $(\mathrm{End}(\mathscr{C}),\circ,\mathrm{Id})$.
\end{example}
As is well-observed in the literature \cite{categories-smc, main-source-js, quantum-ct}, standard monoidal categories satisfy a `strictification theorem' stating that each such category is monoidally equivalent to a strict monoidal category. $\alpha$-ary unbiased monoidal categories follow a similar law; the proofs are very similar (see e.g. Ref. \citenum{quantum-ct}), but it is necessary to define the notions relevant to the theorem.
\begin{definition}[unbiased monoidal functor]
Let $(\mathscr{C}, \mathcal{T}_\bullet, \iota, \eta)$,  $(\mathscr{C}', \mathcal{T}_\bullet', \iota', \eta')$ be unbiased monoidal categories of arity $\alpha$. An \textit{unbiased monoidal functor} between them is a pair $(F, \chi_\bullet)$ consisting of a functor $F\colon\mathscr{C}\to\mathscr{C}'$, and a family of natural isomorphisms
\begin{equation}
\chi_\bullet\colon\alpha\to\bigcup_{\beta<\alpha}{\mathrm{mor}([\mathscr{C}^\beta,\mathscr{C}'])}\quad\text{s.t.}\quad \chi_\beta:\mathcal{T}'_\beta\circ F^\beta\isomto F\circ\mathcal{T}_\beta,
\end{equation}
such that the following diagram commutes for each ordinal partition $p_\bullet\colon\gamma\to\alpha$,
\begin{equation} \label{fnc-cond-1}
\begin{tikzcd}
\mathcal{T}_\gamma'\circ\prod_{i<\gamma}{\left(\mathcal{T}_{p_i}\circ F^{p_i}\right)} \arrow[dd, "\eta'(p_\bullet)F"'] \arrow[rrr, "\mathcal{T}_\gamma' \prod_{i<\gamma}{\chi_i}"] &                         &  & \mathcal{T}_\gamma'\circ F^\gamma\circ\prod_{i<\gamma}{\mathcal{T}_{p_i}} \arrow[dd, "\chi_\gamma \prod_{i<\gamma}{\mathcal{T}_{p_i}}"'] \\
&                         &  &                        \\
\mathcal{T}_\beta'\circ F^\beta \arrow[r, "\chi_\beta"']                                        & F\circ\mathcal{T}_\beta &  & F\circ\mathcal{T}_\gamma\circ\prod_{i<\gamma}{\mathcal{T}_{p_i}} \arrow[ll, "F\eta(p_\bullet)"]    
\end{tikzcd}
\end{equation}
called the \textit{pentagon equation}, as does the following,
\begin{equation} \label{fnc-cond-2}
\begin{tikzcd}
& F \arrow[ldd, "\iota'F"'] \arrow[rdd, "F\iota"] &   \\
&                                                 &   \\
\mathcal{T}_1'\circ F \arrow[rr, "\chi_1"'] &         & F\circ\mathcal{T}
\end{tikzcd}
\end{equation}
called the \textit{triangle equation}.
\end{definition}
As is usually done for monoidal categories, unbiased monoidal functors can be used to define weak and strong equivalences. Similarly to the standard case, we can promote the former to the latter using the axiom of choice \cite{categories-smc}.
\begin{definition}[unbiased monoidal equivalence]
An unbiased monoidal functor $F:\mathscr{C}\to\mathscr{C}'$ is a \textit{weak equivalence} if it is fully faithful and essentially surjective. It is a \textit{strong equivalence} if there is a specified unbiased monoidal functor $G:\mathscr{C}'\to\mathscr{C}$ such that $G\circ F\simeq\mathrm{id}_\mathscr{C}$ and $F\circ G\simeq\mathrm{id}_{\mathscr{C}'}$.
\end{definition}
\begin{proposition}
A weak unbiased monoidal equivalence induces a strong equivalence given the axiom of choice (\texttt{AC}).
\begin{proof}
Let $(F\colon\mathscr{C}\to\mathscr{C}',\chi_\bullet)$ be a unbiased weak monoidal equivalence. For objects $A,B\in\mathscr{C}$, define the equivalence relation $A\sim B\iff FA\simeq FB$. By \texttt{AC} there is a function $f\colon\mathscr{C}/\!\sim\,\to\mathscr{C}$ s.t. $F(f([A]))\simeq FA$; essential surjectivity determines a function $g\colon\mathscr{C}'\to\mathscr{C}/\!\sim$, and thus for each $A'\in\mathscr{C}'$ let $GA'=f(g(A'))$. Let $G$ take morphisms $A'\to B'$ to morphisms $GA'\to GB'$ pre and post-composed by isomorphisms (chosen again by \texttt{AC}). Then $FG$ and $GF$ are both naturally isomorphic to their respective identities (say by maps $\zeta$ and $\xi$, respectively). It remains to show that $G\colon\mathscr{C}'\to\mathscr{C}$ is unbiased monoidal; consider $\forall\beta<\alpha$ the natural isomorphism
\begin{equation}
    G\chi_\beta^{-1}G^\beta\colon GF\circ\mathcal{T}_\beta\circ G^\beta\isomto G\circ\mathcal{T}'_\beta\circ F^\beta G^\beta.
\end{equation}
Noting this we can define a natural isomorphism $\chi'_\beta$ by the commutativity of the following diagram,
\begin{equation}
\begin{tikzcd}
GF\circ\mathcal{T}_\beta\circ G^\beta \arrow[rr, "G\chi_\beta^{-1}G^\beta"] \arrow[dd, "\xi(\mathcal{T}_\beta\circ G^\beta)"] &  & G\circ\mathcal{T}'_\beta\circ F^\beta G^\beta \arrow[dd, "(G\circ\mathcal{T}_\beta')\zeta^\beta"] \\
&  &                                                             \\
\mathcal{T}_\beta\circ G^\beta \arrow[rr, "\chi'_\beta"]         &  & G\circ\mathcal{T}_\beta'                                    
\end{tikzcd}
\end{equation}
which clearly satisfies Eqs. \ref{fnc-cond-1} and \ref{fnc-cond-2} for all $\beta$.
\end{proof}
\end{proposition}
\begin{theorem}[Unbiased Strictification]\label{thm:strictification}
Every $\alpha$-ary unbiased monoidal category is monoidally equivalent to a strict $\alpha$-ary unbiased monoidal category. With \texttt{AC} it is strongly equivalent to one.
\begin{proof}
As discussed, the proof is almost identical the standard case given, e.g., in Ref. \citenum{quantum-ct}. The key observation is that we may recover the unit object, unitors, and associator from the unbiased axioms as follows: the unit is $\mathcal{T}_0(*)$, the left unitor is the composite
\begin{equation*}
\begin{tikzcd}
\mathcal{T}_2\circ(\mathcal{T}_0\times\mathrm{id}) \arrow[rr, "\mathcal{T}_2(\mathrm{id}\times\iota)"] &  & \mathcal{T}_2\circ(\mathcal{T}_0\times\mathcal{T}_1) \arrow[rr, "{\eta_{(0,1)}}"] &  & \mathcal{T}_1 \arrow[rr, "\iota^{-1}"] &  & \mathrm{id},
\end{tikzcd}
\end{equation*}
and the associator is defined in the obvious fashion by noting that each way of ``bracketing'' the two-fold product is naturally isomorphic by $\eta$ to $\mathcal{T}_3$. The right unitor is obtained similarly to the left one with the ordinal partition $(0, 1)$.
\end{proof}
\end{theorem}
Because of this result, studying $\alpha$-ary categories is no different from studying their strict variants (at least regarding \textit{categorical} properties that respect the equivalence principle); we therefore choose the latter as to avoid handling the $\eta$ morphisms.

This framework gives rise to an expected but interesting formulation with infinitary Lawvere Theories (see e.g. Ref. \citenum{lawvere-theories}). Let $\mathbf{Ord}$ be the $\mathbf{Set}$ skeleton of von Neumann ordinals with inclusion functor $\mathcal{I}_\mathbf{Ord}$. Let $\mathcal{F}\colon\mathbf{Mon}\to\mathbf{Set}$ be the forgetful functor; this admits a left adjoint $\mathcal{L}\dashv\mathcal{F}$ sending sets to their corresponding free monoids. Since $\mathcal{L}\circ\mathcal{I}_\mathbf{Ord}$ is injective on objects, its image forms a full subcategory $\mathscr{T}\subseteq\mathbf{Mon}$. Then, ignoring size issues, any functor
\begin{equation}
\mathcal{M}\colon\mathscr{T}^\mathrm{op}\to\mathbf{Cat}
\end{equation}
that preserves small products is identified with a strict $\mathbf{Ord}$-ary unbiased monoidal category. The $\alpha$-ary version can be recovered by replacing $\mathbf{Ord}$ with the subcategory $\alpha\subseteq\mathbf{Set}$ in the previous paragraph.

One can show with some difficulty that each free monoid $M$ in $\mathscr{T}$ can be written distinctly as the sequential limit of a constant functor
\begin{equation}
M\cong\varinjlim{\left(\Delta_M\colon\beta\to\mathscr{T};\gamma\mapsto 1\right)}
\end{equation}
for some $\beta<\alpha$ and the free monoid over a single element $1$, which is just a $\beta$-fold product in $\mathscr{T}^\mathrm{op}$. Then $\mathcal{M}(1)$ is an unbiased $\alpha$-ary monoidal category with product functors
\begin{equation}
\mathcal{T}_\beta=\mathcal{M}\left(1\to\beta;(0)\mapsto (012\dots\beta)\right)
\end{equation}
for all $\beta<\alpha$. As is intuitive, the image of $\mathcal{M}$ is the full subcategory of $\mathbf{Cat}$ that is the aforementioned power category of $\mathcal{M}(1)$.

A basic conclusion from the axioms of a standard monoidal category---yet the one that gives rise to most of its diagrammatic theory---is the \textit{interchange law}, stating that, for properly composable morphisms,
\begin{equation}
(g\circ f)\otimes (j\circ h)=(g\otimes j)\circ (f\otimes h).
\end{equation}
This is a direct consequence of bifunctoriality of $\otimes$. In the infinitary case, however, while we can perform infinitary interchange \textit{internally} to the arguments of the tensor functors, we need additional assumptions to do so on the outside.
\begin{definition}[smooth unbiased category]
An $\alpha$-ary unbiased monoidal category $\mathscr{C}$ is \textit{smooth} if for each object $X$ the functor
\begin{equation}
    \mathcal{T}_2(X,-)\colon \mathscr{C}\to\mathscr{C}
\end{equation}
is sequentially cocontinuous (i.e. preserves sequential colimits). It is \textit{cosmooth} if it preserves cosequentially cocontinuous (i.e. preserves colimits of diagrams $\alpha^\mathrm{op}\to\mathscr{C}$).
\end{definition}
It is not difficult to see that this condition implies that every tensor functor is sequentially cocontinuous in the last argument. If the condition of smoothness holds, we have the following obvious corollary.
\begin{theorem}[Infinitary Interchange Law]
If $\mathscr{C}$ is a smooth $\alpha$-ary monoidal category with products $\mathcal{T}_\bullet$, then the following two conditions hold:
\begin{enumerate}
\item Consider two sequences of morphisms in $\mathscr{C}$, i.e. functors $F,G\colon\beta\to\mathscr{C}$ for the thin category of $\beta<\alpha$ such that both $F,G$ admit transfinite compositions $f\colon F_0\to\varinjlim F$ and $g\colon G_0\to\varinjlim G$ respectively. Then
\begin{equation}\label{inf-int-l1}
\mathcal{T}_2(f,g) = h\colon\mathcal{T}_2(F_0,G_0)\longrightarrow\varinjlim\big(\mathcal{T}_2\circ (F\times G)\big)
\end{equation}
\item If $(a_0,a_1),(b_0,b_1),\dots$ is a morphism sequence of length $\gamma$ such that the composition of each pair exists, then 
\begin{equation}
\mathcal{T}_\gamma\big(a_0\circ a_1,b_0\circ b_1,\dots\big)=
\mathcal{T}_\gamma\big(a_0,a_1,\dots\big)\circ \mathcal{T}_\gamma\big(b_0,b_1,\dots\big)
\end{equation}
\end{enumerate}
\end{theorem}
When such a smooth category is in addition cocomplete, it is a very well-behaved setting for infinitary operations by virtue of the infinite interchange law because all transfinite compositions exist and are compatible with the product. If furthermore the category in question is a saturated subcategory, then all possible infinitary manipulations of bracketings are equalities. This condition, however, is too difficult to satisfy in practice for any category, perhaps other than the thin category $\alpha\subseteq\mathbf{Ord}$ itself.

\subsection{Colimit-expanded Categories}
Beyond the standard examples of categories with all small (co)products, it is rare to naturally encounter infinitary monoidal categories. In addition, the axioms are relatively difficult to show in practice. There is, however, a general process by which one can construct a (smooth) fully infinitary category from a standard one under certain conditions. This construction and its derivatives will not only provide many examples, but also will also give rather profound results about various existing constructions in category theory.

The general idea is to take a monoidal category and use colimits of the finite tensor products to define the functors $\mathcal{T}_\beta$. Although we will not write it out explicitly, whenever we refer to sequential colimits (i.e. those of diagrams $\alpha\to\mathscr{C}$ for an ordinal $\alpha$), we just as well may use \textit{cosequential colimits} taken over diagrams out of $\alpha^{\mathrm{op}}$. 

Now, let $(\mathscr{C},\otimes,I)$ be an ordinary monoidal category. Seemingly the most general case depends on a choice of a cone from $I$ to the identity $\mathrm{Id}\colon\mathscr{C}\to\mathscr{C}$, say $\varepsilon\colon\Delta_I\to\mathrm{Id}$. We proceed with transfinite recursion. Suppose we have a sequence $(X_\gamma)_{\gamma<\alpha}\in\mathscr{C}^\alpha$; to define the infinitary tensor product $\mathcal{T}_\alpha$, we first recursively let $F(0)=I$ and $F(\gamma+1)=F(\gamma)\otimes X_\gamma$. Under $F$ let successor maps $\gamma\to\gamma+1$ go to the composite
\begin{equation}
    F(\gamma)\xrightarrow[\hspace*{25px}]{\rho^{-1}}F(\gamma)\otimes I\xrightarrow[\hspace*{25px}]{1\otimes \varepsilon}F(\gamma)\otimes X_\gamma,
\end{equation}
and finally extend $F$ to a functor $\alpha\to\mathscr{C}$ by taking
\begin{equation}\label{eq:exp-by-col}
    F(\beta)=\varinjlim_{\gamma<\beta}F_\gamma
\end{equation}
for any limit ordinal $\beta$. Because of the associator and unitors of the underlying category $\mathscr{C}$, the axioms of a fully infinitary category are almost satisfied, save for the fact that the coherence conditions cannot be proven in the limit ordinal case. To resolve this, we must impose the additional requirement that for any $F\colon\alpha\to\mathscr{C}$, there is an is an isomomorphism
\begin{equation}
X\otimes\varinjlim_{\gamma<\alpha}{F(\gamma)}\cong\varinjlim_{\gamma<\alpha}\big(X\otimes F(\gamma)\big),
\end{equation}
i.e. $\mathscr{C}$ is smooth. Moreover, smoothness guarantees that each colimit of the form in Eq. \ref{eq:exp-by-col} exists. Notably, the above construction is just an example of transfinite compositions of the maps $\varepsilon$, which serve to `glue' the full product together by partial products. This gives rise to the following definition.
\begin{definition}[colimit expansion]
Let $\mathscr{C}$ be a monoidal category equipped with a cone $\varepsilon$ from $I$ to $\mathrm{Id}$. $\mathscr{C}$ is  \textit{colimit-expandable} if for each object $X$ of $\mathscr{C}$ the functor $X\otimes -\colon \mathscr{C}\to\mathscr{C}$ is sequentially cocontinuous. The resulting smooth, fully infinitary category as constructed above is the \textit{colimit expansion} of $\mathscr{C}$, and the pair $(\mathscr{C}, \varepsilon)$ is termed a \textit{colimit-expanded} category.
\end{definition}
The obvious examples are $\mathbf{Set}$, $\mathbf{Rel},\mathbf{Cat}$, etc. all equipped with their coproducts and initial objects as units. Recall that we just as easily may use the product if instead a cone from $\mathrm{Id}$ to $I$ is employed. $\mathbf{Top}$ is an interesting non-example: sums of topological spaces do not extend canonically to the infinitary case. As we will see later, however, the ``convenient category of topological spaces" \textit{is} (dually) colimit-expandable.
\begin{proposition}
A colimit-expanded category $(\mathscr{C}, \varepsilon)$ is necessarily semicocartesian (i.e. $I$ is initial). Thus the legs of $\varepsilon$ are trivially the unique maps out of $I$.
\begin{proof}
First notice that from the definition of a fully infinitary category, in particular Eq. \ref{eq:unb-ax-2}, it follows that for any $\alpha$ there is an isomorphism $\mathcal{T}_\alpha(I,I,\dots)\cong I$ induced by $\eta$ and $\iota$. Now suppose that the leg $\varepsilon_I\colon I\to I$ was not an isomorphism and let $\alpha$ be a limit ordinal; since the product functor $\mathcal{T}_\alpha$ is constructed via transfinite compositions of $\epsilon$ which reflect isomorphisms, we must have $I\ncong\mathcal{T}_\alpha(I,I,\dots)$, a contradiction.

Then since $\varepsilon_I\colon I\cong I$, the commutativity condition of $\varepsilon$ guarantees that any two morphisms $f,g\colon I\to X$ are equal. Of course each $I\to X$ exists by definition of $\varepsilon$, so $I$ is initial in $\mathscr{C}$.
\end{proof}
\end{proposition}
Because of this, there is not, in fact, freedom to choose various cones $I\to\mathrm{Id}$, and thus we need not mention $\varepsilon$ hereafter. The following statement, though a direct consequence of the limit preservation properties of adjunctions, creates many of the rest of the examples of smooth monoidal categories. We recall that a closed monoidal category is a symmetric one in which, for each object $X$, the functor $-\otimes X$ has a right adjoint $[X,-]\colon\mathscr{C}\to\mathscr{C}$ called the `internal hom'.
\begin{theorem}
Any \textit{semicocartesian closed} category is colimit-expandable. In particular, any cocartesian closed category is colimit-expandable. By the same argument with cosequential colimits, semicartesian closed categories, and thus cartesian closed categories, are as well colimit-expandable.
\end{theorem}
Firstly, this gives an interesting result about closed monoidal categories; it implies a rather unintuitive underlying relationship between the existence of an internal hom and the possibility and uniqueness of an infinitary extension. These semicocartesian closed categories have various special properties that might give rise to such a relationship (e.g., by the adjunction formula, $T\cong[X,I]$); the author believes that a study of these categories in their own right would be worthwhile.

Secondly, it provides us with many obvious examples, circumventing the tedious axioms. For instance, we can remedy the earlier non-example of $\mathbf{Top}$ by taking the cartesian closed subcategory $\mathbf{CGHaus}$ of compactly-generated Hausdorff spaces \cite{categories-smc}. In addition, we immediately have that all topoi are colimit-expandable and thus uniquely smooth infinitary categories, giving yet another useful characteristic of topoi that may have future application.

There are some examples of semi(co)cartesian closed categories that are not cartesian closed; although most of these perhaps are slightly artificial, such examples show that the colimit-expansion is versatile, and need not deal only with the categorical product. The primary example, thought, that the author has found in usage is that of nominal sets (see Andrew Pitts' \textit{Nominal Sets} \cite{nominal-sets}; we follow his work here).
\begin{example}
For any set $A$, by $\mathrm{Perm}(A)$ we denote the subgroup of $\mathrm{Sym}(A)$ containing only those permutations that leave only finitely many elements of $A$ changed. Fix a countably infinite $\mathbb{A}$, and let $X$ be a $\mathrm{Perm}(\mathbb{A})$-set (i.e. an object in the $\mathrm{Perm}(\mathbb{A})\mathbf{Set}$ topos). We say that a subset $A\subseteq\mathbb{A}$ \textit{supports} an $x\in X$ if $\forall\pi\in\mathrm{Perm}(\mathbb{A})$, we have $\left. \pi\right\vert_A=1_A\implies \pi\cdot x$. Any $x$ is \textit{finitely supported} if there is a finite subset of $\mathbb{A}$ which supports it. If it is so, then we write the intersection of all finite supports as $\mathrm{supp}_X{x}$ Then $\mathbf{Nom}$ is defined as the full subcategory of $\mathrm{Perm}(\mathbb{A})\mathbf{Set}$ consisting of $\mathrm{Perm}(\mathbb{A})$-sets in which each element is finitely supported (Ref. \citenum{nominal-sets}, Ch. 2.1).

$\mathbf{Nom}$ inherits a cartesian closed structure from the topos $\mathrm{Perm}(\mathbb{A})$ with the categorical product, so it is colimit-expandable under thereunder. However it is also as such under the \textit{separated product} (Ref. \citenum{nominal-sets}, Ch. 3.4): we define the `freshness relation' $\#_{X,Y}\subseteq X\times Y$ by $x\,\#_{X,Y}\, y \iff \mathrm{supp}_X{x}\,\cap\,\mathrm{supp}_Y{y}=\varnothing$. The relation is itself a nominal set, giving a bifunctor $(X,Y)\mapsto\,\#_{X,Y}$ called the separated product, under which $\mathbf{Nom}$ is made monoidal with the terminal object as the unit. Finally by Ref. \citenum{nominal-sets}, Thm. 3.12, it is semicartesian closed.
\end{example}
We conclude with a general demonstration of obtaining semi(co)carte- sian closed categories, and thus colimit-expandable ones, from monoidally closed ones via the (co)slice category. 
\begin{proposition}
Let $\mathscr{C}$ be closed closed monoidal with unit $I$. Both the slice category $\mathscr{C}\slash I$ and the coslice category $I\slash \mathscr{C}$ are colimit-expandable.
\begin{proof}
It is sufficient to consider the slice case. Observe that there is a natural monoidal structure inherited from $\mathscr{C}$ on the slice category given by
\begin{equation}
(X,f)\otimes (Y,g)=(X\otimes Y,f\otimes g\colon X\otimes Y\to I\otimes I\cong I),
\end{equation}
in which $(I, 1_I)$ is terminal, implying that $\mathscr{C}\slash I$ is semicartesian. Next, since $\mathscr{C}$ is closed, it has a unique internal hom
\begin{equation}
[-,-]\colon\mathscr{C}^{\mathrm{op}}\times\mathscr{C}\to\mathscr{C},
\end{equation}
so we may define a similar bifunctor for the slice category by
\begin{equation}
(X,f),(Y,g)\mapsto \big([X,Y], [f,g]\big).
\end{equation}
We wish to show that the pair
\begin{equation}
-\otimes\,(X,f)\dashv\big[(X,f),-\big]
\end{equation}
is an adjunction for each $(X,f)$ in $\mathscr{C}\slash I$. Quite simply, given a morphism of $\mathscr{C}\slash I$ from $(X\otimes Y,f\otimes g)$ to some $(Z,h)$, we get a unique corresponding map given by the adjunct. This is a valid morphism in the slice category as it makes the diagram
\begin{equation}
\begin{tikzcd}
X \arrow[d, "f"'] \arrow[r, "u"] & {[Y,Z]} \arrow[d, "{[g,h]}"] \\
I \arrow[r, "\lambda^{-1}"']     & I\otimes I         
\end{tikzcd}
\end{equation}
commute by the naturality of the relevant adjunction in $\mathscr{C}$, which also gives the naturality of this bijection.
\end{proof}
\end{proposition}
It is interesting that even though a general closed monoidal category is not colimit-expandable, the corresponding (co)slice is. This suggests that `in spirit' all categories that are closed monoidal have the \textit{structure} to be infinitary monoidal, but rather it is a technicality of our definition that they are not so. Future work, perhaps, could amend the axioms to resolve this.

\section{Monoidal Diagrams}\label{monoidal-diagrams}
As discussed, the use of diagrams similar to those employed in monoidal category theory dates back to the work of Penrose \cite{penrose-1,penrose-2}. When applied specifically to monoidal categories, they become very useful because isotopic deformations preserve represented morphisms. 

In general, each diagram represents a morphism constructed by compositions and tensor products of other morphisms; these morphisms are the nodes of the diagram, while the edges are the objects. Composition takes place vertically while the product acts horizontally. Constructing the diagram for two morphisms and seeing if they are isotopic to one another provides an efficient way to determine the equivalence of two morphisms, which is often non-trivial given the monoidal product \cite{survey-graph-types}.

In reality, the use of these `string diagrams' in category theory is not restricted to monoidal categories. In general, any weak 2-category operates in a manner respectful of string diagrams \cite{string-diagrams} with 0-cells represented by areas, 1-cells by lines, and 2-cells by points. Since monoidal categories are weak 2-categories with one object, the formalism simplifies considerably.

These string diagrams are usually topological entities which map to categories \cite{main-source-js}. While this framework is certainly the most natural and simple to prove soundness/correctness within, it makes it difficult to generalize to the strictly categorical notions we wish to consider, i.e. transfinite composition and infinitary monoidal product. To remedy this, we shall use define a diagram purely algebraically. This not only strengthens its categorical properties, but permits for its use as a mathematical object in its own right outside of calculational purposes.

The next few sections will look at, in detail, our upgraded definition for monoidal diagrams. Sec. \ref{mon-di-def} focuses on defining diagrams on a monoidal category, exemplifying this concept, and proving some preliminary results therefor. Secs. \ref{mon-di-unr} and \ref{mon-di-rea} will both be devoted to understanding the morphism which a given diagram represents; the ultimate goal will be to define a category of diagrams, and then establish a functorial correspondence between it and the underlying category.

\subsection{Definitions} \label{mon-di-def}
We shall construct monoidal diagrams as algebraic structures equipped with two relations; we begin with the necessary definitions.
\begin{definition}[diagrammatic relation]
A homogeneous relation $R\subseteq X\times X$ is diagrammatic if the following conditions hold $\forall x,y\in X$:
\begin{enumerate}
    \item \textit{(weak totality)} $\exists z\in X\,\big(xRz\lor zRx\big)$
    \item \textit{(antisymmetry)} $xRy\implies\lnot yRx$
    \item \textit{(irreflexivity)} $\lnot xRx$
\end{enumerate}
\end{definition}
Note that the last two conditions could be replaced with that of asymmetry; a diagrammatic relation is equivalently an asymmetric relation where each element is comparable to another.
\begin{definition}[induced first-order relation] We call a homogeneous relation $R$ on a set $X$ a \textit{first-order relation}, a homogeneous relation $S\subseteq R\times R$ on $R$ a \textit{second-order relation}, ad infinitum. For a second-order relation, we define the \textit{induced first-order relation} of $S$ as
\begin{equation}
    S^*=\bigcup_{i=1,2}{\pi^2_i[S]},
\end{equation}
where $\pi_i$ is the $i$th projection map of the Cartesian product $R\times R$, $f^2\colon A\times A\to B\times B$ represents the Cartesian square of a function, and the notation $f[X]$ for a function $f\colon A\to B$ refers to the set $\{f(x):x\in X\subseteq A\}$. 
\end{definition}
Intuitively, if the relation $S$ in the above orders the edges of a digraph, then $S^*$ orders the nodes based on their membership of an edge; this is the case in which we will employ the definition.
\begin{definition}[conditional construction]\label{rec-constr}
Suppose that a second order-relation $S$ on a set $X$ is a homogeneous relation on $R$. Then $S$ satisfies the axiom of \textit{conditional construction} if
\begin{equation}
    \forall P\in R^2\,\left(\bigvee_{i,j\leq 2}{(\pi_i\times\pi_j)(P)\in S^*}\implies P\in S\right)
\end{equation}
holds; in this case, $S$ is \textit{conditionally constructed}.
\end{definition}
We will use the axiom of conditional construction simply to ensure the coherence of a relation used to order the edges of a digraph as discussed above. With this setup, we can introduce monoidal diagrams. For what follows, fix a strict unbiased monoidal category $\mathscr{C}$.
\begin{definition}[monoidal diagram]\label{mon-dia}
A (monoidal $\mathscr{C}$-diagram $(\mathcal{D},E,H,\mu)$ consists of a set $\mathcal{D}$, a diagrammatic relation $E$ on $\mathcal{D}$ called the \textit{vertical comparator}, a partial ordering $H$ on $E$ called the \textit{horizontal comparator} (or edge set), and a function $\mu\colon\mathcal{D}\to\mathrm{mor}(\mathscr{C})$. These data are required to satisfy the following conditions:
\begin{enumerate}
    \item $E, H$ are well-founded relations
    \item Each $E$-triangle (set of two $E$-elements which take the same value on $\pi_1$ or $\pi_2$) is $H$-comparable
    \item Any two $E$-minimal elements are $H^\star$-comparable 
    \item $H$ is conditionally constructed
\end{enumerate}
\end{definition}
There is much to be said about this definition. The four conditions are non-trivial, and it is not immediately clear that this definition, taken as an axiomatic system, is consistent, complete, or independent, the first two being quite important in this context. As we develop the theory further, however, the three conditions will hold apparently.

Separately, Def. \ref{mon-dia} is not entirely useful until we prescribe a method for actually `viewing' such a diagram, in the literal sense. To do so one begins by laying out the nodes in $\mathcal{D}$ vertically by $E$, and then connecting with a segment any $E$-related nodes. Next, these segments are ordered horizontally by $H$; it is not difficult to see why this algebraic structure is resistant to isotopic deformations.
\begin{definition}[isomorphism of monoidal diagrams]
Let $(\mathcal{D},E,H,\mu)$, $(\mathcal{D}',E',H',\mu')$ be monoidal $\mathscr{C}$-diagrams. An \textit{isomorphism of diagrams} is a bijection $\delta\colon\mathcal{D}\to\mathcal{D}'$ such that 
\begin{equation}
    \delta\colon(\mathcal{D},E)\to(\mathcal{D}',E')
\end{equation}
is a digraph isomorphism, the induced bijection $\delta'\colon E\to E'$ makes
\begin{equation}
    \delta'\colon (E,H)\to (E',H')
\end{equation}
an order isomorphism, and finally $\mu'\circ\delta=\mu$. These conditions can be summarized by the commutative diagram below.
\begin{equation*}
\begin{tikzcd}
H \arrow[dd, dashed] \arrow[r, hook] & {\mathcal{D}^2}^2 \arrow[r, harpoon, bend right] \arrow[r, harpoon, bend left]  & E \arrow[dd, dashed] \arrow[r, hook] & \mathcal{D}^2 \arrow[r, harpoon, bend left] \arrow[r, harpoon, bend right]  & \mathcal{D} \arrow[dd, "\delta"] \arrow[rrd, "\mu"] &  &                           \\
&                                                                     &                                      &                               &                                                     &  & \mathrm{mor}(\mathscr{C}) \\
H' \arrow[r, hook]                   & {\mathcal{D}'^2}^2 \arrow[r, harpoon, bend left] \arrow[r, harpoon, bend right] & E' \arrow[r, hook]                   & \mathcal{D}'^2 \arrow[r, harpoon, bend left] \arrow[r, harpoon, bend right] & \mathcal{D}' \arrow[rru, "\mu'"]                    &  &                          
\end{tikzcd}
\end{equation*}
Here, the arrows $\rightharpoonup$, $\hookrightarrow$, $\dashrightarrow$ represent projections, inclusions, and induced bijections respectively.
\end{definition}
We acknowledge that isomorphism of diagrams is an extremely strong condition: diagrams are isomorphic if and only if they are \textit{exactly} the same aside from their nodes being labeled differently.

Diagrams may be adjoined to one another using the coproduct on the underlying set and relations and the product functors in the underlying monoidal category. Specifically, for an $\alpha$-sequence $\{(\mathcal{D}_\gamma,E_\gamma,H_\gamma,\mu\gamma):\gamma<\alpha\}$, we set
\begin{equation}
\coprod_{\gamma<\alpha}{(\mathcal{D}_\gamma,E_\gamma,H_\gamma,\mu_\gamma)}=\left(\bigsqcup_{\gamma<\alpha}\mathcal{D}_\gamma,\coprod_{\gamma<\alpha}E_\gamma,\coprod_{\gamma<\alpha}H_\gamma,{\langle{\mu_\gamma}\rangle}_{\gamma<\alpha}\right),
\end{equation}
where the fourth component on the right-hand side is the unique map induced by the coproduct. We call this process \textit{attachment}. It is worth noting that the attachment of diagrams coincides with the true coproduct in the category
\begin{equation}
\mathbf{Set}\times\mathbf{Rel}\times\mathbf{Pos}\times(\mathbf{Set}\slash\mathrm{mor}(\mathscr{C})).
\end{equation}
This construction will be useful later in setting up a category of all diagrams. Hereafter, we will turn our attention to understanding the internal structure of monoidal diagrams. For what follows, fix a $\mathscr{C}$-monoidal diagram $(\mathcal{D}, E, H, \mu)$.
\begin{proposition}
Given a well-founded diagrammatic relation $E$ on $D$, the reflexive transitive closure $E^+_r$ is a well-founded partial order.
\begin{proof}
$E^+_r$ is trivially a partial order since $E$ is already antisymmetric. Let $A\subseteq D$, and let $a\in A$ be a $E$-minimal element (which exists by well-foundedness). We wish to show that $a$ is $E^+_r$-minimal (and thus $E^+_r$ is well-founded); if it weren't, then $\exists b\neq a\in A$ s.t. $bE^+_r a$. Since $\lnot bEa$, $\exists c\in A$ s.t. $bEc$ and $cEa$, but then $a$ is not $E$-minimal, a contradiction.
\end{proof}
\end{proposition}

\begin{definition}[ordinal rank of a diagram]
Let $E^+_r$ be the reflexive transitive closure of $E$, and $\mathfrak{C}_m(E^+_r)$ be the set of $\subseteq$-maximal chains therein. Each $C\in\mathfrak{C}_m(E^+_r)$ is well-ordered, and thus isomorphic to a unique ordinal, its ordinal rank, denoted $\mathrm{rank}(C)$. The rank of $\mathcal{D}$ is then the largest such rank:
\begin{equation}
    \mathrm{rank}(\mathcal{D})=\max{\{\mathrm{rank}(C):C\in\mathfrak{C}_m(E^+_r)\}}.
\end{equation}
\end{definition}
\begin{definition}[segmentation, layering]
Let $\alpha=\mathrm{rank}(\mathcal{D})$. The \textit{segmentation} of $\mathcal{D}$ is the function
\begin{equation}
    S_\bullet\colon\alpha\to\mathcal{P}(\mathcal{D})
\end{equation}
defined by transfinite recursion as follows: $S_0=\varnothing$, and $\forall\beta<\alpha$
\begin{equation}
    S_\beta=\big\{m\in\mathcal{D}:m\text{ is minimal in }\mathcal{D}\setminus\Big(\bigcup_{\gamma<\beta}{S_\gamma}
    \Big)\big\},
\end{equation}
which always exists by well-foundedness. We define the closely related concept of the \textit{layering} of $\mathcal{D}$, which is the map
\begin{equation}
\begin{split}
    \mathcal{S}\colon\alpha&\to\mathcal{P}(\mathcal{D}). \\
                              \beta&\mapsto\bigcup_{\gamma<\beta}{S_\gamma}
\end{split}    
\end{equation}
Then we may more simply write $S_\beta$ as the set of minimal elements of the set difference $\mathcal{D}\setminus\mathcal{S}(\beta)$.
\end{definition}
\begin{proposition}
The segmentation of $\mathcal{D}$ partitions its nodes.
\begin{proof}
Suppose $S_\bullet$ wasn't a disjoint family; then $\exists \beta\neq\gamma$ such that $S_\beta\cap S_\gamma\neq\varnothing$, and $\exists x\in S_\beta\cap S_\gamma$. Without loss of generality, assume $\gamma<\beta$; $x\in S_\gamma\implies x\in \mathcal{S}(\beta)$, which contradicts $x\in S_\beta$. We next must show that the extension $\mathcal{S}(\alpha+1)=\mathcal{D}$. If this was false, then for $m\in S_\beta$, $m$ would not be a minimal element, giving the required contradiction.
\end{proof}
\end{proposition}

\begin{remark}\label{lone-node-dilemma}
At present, `lone' nodes (those not connected by an edge to any other node) are \textit{not} permissible in these diagrams, as they break the weak totality condition on $E$. We do not consider this a large problem for multiple reasons. Firstly, they could relatively easily be integrated into the theory, however the formalism would lose a certain elegance it has, that of the actual tuples storing only necessary information to reconstruct the original morphisms. More importantly, this condition is imposed without loss of any generality because, given a diagram that does \textit{not} satisfy this weak totality condition, we could add a layer below all of the others and affix identity morphisms to the lone nodes, promoting it to a permissible monoidal diagram.
\end{remark}

\subsection{Unresolved Edges} \label{mon-di-unr}
The idea of `unresolved edges' presents an initial problem in this diagrammatic theory. When defining how to \textit{read} a diagram, we shall break it down into its segmentation. But, for morphisms to compose properly, we must insert an identity edge to make uniform the width of each each layer. The following construction will fix this issue; once again, fix a monoidal diagram $(\mathcal{D},E,H,\mu)$.
\begin{definition}[$\gamma$-unresolved edge] An edge $e\in E$ is \textit{$\gamma$-unresolved} if $\pi_1(e)\in\mathcal{S(\beta)}$ and $\pi_2(e)\in\mathcal{D}\setminus\mathcal{S}(\beta)$.    
\end{definition}
\begin{definition}[$e$-incision, cohesion map] \label{inc-coh}
Suppose $\mathcal{D}$ has an unresolved edge $e$ (at any $\gamma$). Let $\mathcal{I}=\{(\pi_1(e), 1_e), (1_e, \pi_2, e)\}$ be an \textit{$e$-incision}, and $E'=E\setminus\{e\}\cup\mathcal{I}$ be the \textit{$e$-incision of $E$}. We then define the \textit{cohesion map} $c\colon E'\to E$ as the function such that $c\vert_{E'\setminus\mathcal{I}}$ is the identity, and $c[\mathcal{I}]=\{e\}$.
\end{definition}
\begin{proposition} \label{unq-ext}
Given the cohesion map $c\colon E'\to E$ as in Def. \ref{inc-coh}, there is a unique partial ordering $H'$ such that the extension $c\colon (E',H_{st}')\to (E, H_{st})$ is order-reflecting (where $A_{st}$ is the corresponding strict partial order of $A$).
\begin{proof}
This immediately follows from the definition of an order-reflecting map. Uniqueness is also preserved by the strictness, since the map $A\mapsto A_{st}$ is bijective on partial orders.
\end{proof}
\end{proposition}
\begin{proposition} \label{coh-diag}
Proceeding with the language of Def. \ref{inc-coh} and Prop. \ref{unq-ext}, let $\mathcal{D}'=\mathcal{D}\cup\{1_e\}$ and $\mu':\mathcal{D}'\to\mathrm{mor}(\mathscr{C})$ such that $\mu'\vert_{\mathcal{D}}=\mu$ and $\mu'(1_e)=(\mathrm{cod}\circ\mu\circ\pi_1)(e)$. Then $(\mathcal{D}',E',H',\mu')$ forms a monoidal diagram.
\begin{proof}
This is an obvious consequence of the fact that $(\mathcal{D},E,H,\mu)$ is itself a monoidal diagram.
\end{proof}
\end{proposition}
\begin{definition}[$e$-resolution, local resolver]
Let $\mathfrak{D}$ be the class of all $\mathscr{C}$-diagrams, and for each diagram $\mathcal{D}$ let $\mathcal{U}^\gamma_\mathcal{D}$ be the set of $\gamma$-unresolved edges in $\mathcal{D}$, and let $\mathcal{U}_\mathcal{D}=\bigcup_{\gamma<\mathrm{rank}(\mathcal{D})}{\mathcal{U}^\gamma_\mathcal{D}}$ be the set of \textit{all} unresolved edges in $\mathcal{D}$. Following Prop. \ref{coh-diag}, we call the diagram $\mathrm{res}(\mathcal{D},e)=(\mathcal{D}',E',H',\mu')$ the \textit{$e$-resolution} of $\mathcal{D}$. We thus have a well-defined class function
\begin{equation}
\mathrm{res}\colon\bigsqcup_{\mathcal{D}\in\mathfrak{D}}{\mathcal{U}_\mathcal{D}}\to\mathfrak{D}
\end{equation}
called the \textit{local resolver} of $\mathscr{C}$-diagrams.
\end{definition}
\begin{proposition}
For a monoidal diagram $\mathcal{D}$, the set $\mathcal{U}_\mathcal{D}$ can be well-ordered without invoking \texttt{AC}.
\begin{proof}
We construct a well-ordering on $\mathcal{U}_\mathcal{D}$ as follows: take the lexicographic ordering $E\times E$ on $\mathcal{D}^2$; since $\mathcal{U}_\mathcal{D}\subseteq E\subseteq\mathcal{D}^2$, there is a well-founded poset uniquely defined through $E\times E$ which exists and satisfies the preceding condition on $\mathcal{U}_\mathcal{D}$.
\end{proof}
\end{proposition}
Because of the existence of this well-ordering, we can well-define the following piecewise function:
\begin{equation}
\begin{split}
\mathrm{res}_*\colon\mathfrak{D}&\to\mathfrak{D} \\
\mathcal{D}&\mapsto
\begin{cases} 
      \mathrm{res}(\mathcal{D},\min{\left(\mathcal{U}_\mathcal{D}\right)}) & \mathcal{U}_\mathcal{D}\neq\varnothing \\
      \mathcal{D} & \mathcal{U}_\mathcal{D}=\varnothing \\
\end{cases}
\end{split}
\end{equation}
This map allows us to, edge-by-edge, resolve all unresolved edges, until it becomes idempotent (i.e. the entire diagram has been resolved).
\begin{definition}[resistivity, global resolver, resolution]
The \textit{resistivity}, $\mathfrak{r}\colon\mathfrak{D}\to\mathbf{Ord}$, of a diagram $\mathfrak{r}(\mathcal{D})$ is the smallest ordinal $\alpha$ such that
\begin{equation}
\mathrm{res}^{(\alpha+1)}_*(\mathcal{D})=\mathrm{res}^{(\alpha)}_*(\mathcal{D}).
\end{equation}
We define the \textit{global resolver} as
\begin{equation}
\begin{split}
\mathrm{res}\colon\mathfrak{D}&\to\mathfrak{D} \\
\mathcal{D}&\mapsto\mathrm{res}_*^{(\mathfrak{r}(\mathcal{D}))}(\mathcal{D}).
\end{split}
\end{equation}
We call $\mathrm{res}(\mathcal{D})$ the \textit{resolution} of the diagram $\mathcal{D}$, and define the subclass $\overline{\mathfrak{D}}$ as the image of $\mathrm{res}$, elements of which are resolved diagrams.
\end{definition}
Having resolved these edges, each diagram can be deconstructed into layers properly; this will be demonstrated in the following section.

\subsection{Readouts \& Validity} \label{mon-di-rea}
We may move towards extracting the morphism encoded by a diagram, which will be termed its \textit{readout}. Before doing so, however, there is an important subtlety to note: as we have not strongly restricted the map $\mu$, it is possible that the morphisms dictated to compose by the structure of a diagram are simply unable to do so due to their domains and codomains; we would like to think of such diagrams as \textit{invalid}. To rigorously address this, we will first determine the potential readout of a diagram, and then formalize its validity. To do so, the following result will be useful.
\begin{proposition}\label{well-ord-seg}
Let $(\mathcal{D},E,H,\mu)$ be a monoidal diagram such that $\mathcal{D}\in\bar{\mathfrak{D}}$. For any $\beta<\mathrm{rank}(\mathcal{D})$, let
\begin{equation}
    H^*_\beta = H^*\cap\big( S_\beta\times S_\beta\big).
\end{equation}
Then $(S_\beta, H^*_\beta)$ is a well-ordered set.
\begin{proof}
Suppose $H^*$ was not well-founded; then there is a non-empty subset $U\subseteq\mathcal{D}$ with no minimal element. Let $U'=\{e\in E:\bigvee_{i=1,2}{\pi_i(e)\in U}\}$. By condition 1 of monoidal diagrams, $(E,H)$ is a well-founded set, and thus $U'$ has a minimal element, say $(c,d)$. Choose $x\in \{c,d\}\cap U$; since $U$ has no minimal element, $\exists$ distinct $y\in U$ such that $y H^* x$. By weak totality, $y$ must be part of at least one pair in $U'$, say $(a,b)\in U'$ that is unique from $(c,d)$. But since $H$ is well-constructed as a second-order relation, this implies that $(a,b)H(c,d)$, which contradicts $(c,d)$ being a minimal element, showing the well-foundedness of $H^*$.

Since $H^*_\beta\subseteq H^*$ and is therefore well-founded itself, it remains to show that $S_\beta$ is totally ordered. We will proceed with transfinite induction on $\beta$; thus first assume $(S_\beta, H_\beta^*)$ is totally ordered. For the successor ordinal case, let $x,y\in S_{\beta+1}$. Note that there must exist $a\in S_\beta$ such that $a E x$; if not, then $x$ would have been minimal in $\mathcal{D}\setminus\mathcal{S}(\beta)$ and thus a member of $S_\beta$ which it cannot be since $S_\bullet$ partitions $\mathcal{D}$. Similarly, there is a $b\in S_\beta$ such that $b E y$. By supposition either $aH^*_\beta b$ or $bH^*_\beta a$, or identically $aH^* b$ or $bH^* a$. But well-constructedness then implies either that $(a,x)H(b,y)$ or $(b,y)H(a,x)$, which by definition means that $x$ and $y$ are $H^*$ comparable and thus obviously $H^*_{\beta+1}$ comparable.

For the limit ordinal case, suppose $x,y\in\beta$ for some limit ordinal $\beta$, and that $(S_\delta,H^*_\delta)$ is totally ordered for all $\delta<\beta$. By the same argument as above, $\exists a,b\in\mathcal{S}(\beta)$ such that $aEx$ and $bEy$. Furthermore, $a\in S_{\delta_1}$ and $b\in S_{\delta_2}$ for $\delta_1,\delta_2<\beta$; if $\delta_1\neq\delta_2$, say without loss of generality that $\delta_1>\delta_2$, then the diagram would contain a $\delta_1$-unresolved edge, which is a contradiction since $\mathcal{D}\in\bar{\mathfrak{D}}$. Thus $\delta_1=\delta_2=\delta$, and so $x$ and $y$ are $H^*_{\beta}$ comparable by the same logic as above. Finally, the base case $\beta=1$ is trivial since $S_1$ consists of all $E$-minimal elements, but by condition three of a monoidal diagram, these must be $H^*$ comparable (and thus $H^*_\beta$ comparable as well). $S_\beta$ is therefore well-ordered.
\end{proof}
\end{proposition}
\begin{remark}
Actually, for the above result, we need not restrict ourselves to the subclass $\bar{\mathfrak{D}}$ of resolved diagrams; doing so does greatly simplify the proof however, especially since we will only need this well-ordering for resolved diagrams. The more general form can be proved by noting that any two nodes in a segment either trace back to an $E$-triangle at some point in the diagram by chain, or they lead to $E$-minimal elements; thus their comparability recurses back up through the diagram by means of the axiom of conditional construction.
\end{remark} 
As discussed, the segments of this diagram indicate where to apply composition and monoidal product. The segmentation is already constructed to be well-ordered, in part, by the vertical comparator. Given now that the horizontal comparator well-orders each segment, we are aware in what order to carry out monoidal product when constructing the readout.
\begin{definition}[reading]
Let $(\mathcal{D},E,H,\mu)$ be a monoidal diagram such that $\mathcal{D}\in\bar{\mathfrak{D}}$. Let $\ell_{\beta}=\mathrm{rank}(S_\beta)$ for all $\beta$. Let
\begin{equation}
S^\bullet_\beta\colon\ell_\beta\isoto S_\beta
\end{equation}
be the unique order-isomorphism between the well-ordering and its ordinal rank, by means of Prop. \ref{well-ord-seg}. Using the notation of Sec. \ref{sec-def-unb}, the \textit{reading} of $\mathcal{D}$ is the ordinal sequence
\begin{equation}
\begin{split}
\mathrm{read}(\mathcal{D})_\bullet\colon\mathrm{rank}(\mathcal{D})&\to\mathrm{mor}(\mathscr{C})\\
\beta &\mapsto \mathcal{T}_{\ell_\beta}\circ\prod_{\gamma<\ell_\beta}{\mu\big(S_\beta^\gamma\big)}.
\end{split}
\end{equation}
\end{definition}
Since we now have a sequence of morphisms in the category which are to be composed, we must impose the validity requirement before doing so.
\begin{definition}[valid diagram]
Let $(\mathcal{D},E,H,\mu)$ be a monoidal diagram such that $\mathcal{D}\in\bar{\mathfrak{D}}$. Then let
\begin{equation}
\mathrm{inv}(\mathcal{D}) = \{\gamma<\mathrm{rank}(\mathcal{D}):\mathrm{dom}(\mathrm{read}(\mathcal{D})_{\gamma+1})\neq \mathrm{cod}(\mathrm{read}(\mathcal{D})_{\gamma})\}
\end{equation}
We shall write $\mathrm{val}(\mathfrak{D})$ for the preimage of $\varnothing$ under $\mathrm{inv}\circ\mathrm{res}\colon \mathfrak{D}\to\mathcal{P}(\mathbf{Ord})$. An element of $\mathrm{val}(\mathfrak{D})$ is a \textit{valid} diagram.
\end{definition}
It is not in general true that $\mathrm{val}(\mathfrak{D})\subseteq\bar{\mathfrak{D}}$; the $\mathrm{res}$ map allows us to consider diagrams valid even if they contain unresolved edges. With this, we intuitively would like to define the readout map as the transfinite composition (as in Def. \ref{def:transfinite-composition}) of the reading. To do so, consider for the diagram $\mathcal{D}\in\mathrm{val}(\mathfrak{D})$ a functor
\begin{equation}
\mathcal{D}_\bullet\colon\mathrm{rank}({\mathcal{D}})\to\mathscr{C},
\end{equation}
such that $\mathcal{D}_\bullet = \mathrm{dom}\circ\mathrm{read}(\mathcal{D})_\bullet$ on objects. On the morphism $\beta\to\beta+1$, let $\mathcal{D}_\bullet$ take the value $\mathrm{read}(\mathcal{D})_\beta$. The composition law is satisfied since $\mathcal{D}$ is a valid diagram.
\begin{definition}[readout]
The \textit{readout} is the function $\mathcal{R}\colon\mathrm{val}(\mathfrak{D})\to\mathrm{mor}(\mathcal{C})$ which maps a diagram $\mathcal{D}$ to the unique morphism induced by the colimit:
\begin{equation}
\begin{tikzcd}
\mathcal{D}_0 \arrow[rr, "\mathcal{R}(\mathcal{D})", dashed] &  & \varinjlim{\mathcal{D}_\bullet}
\end{tikzcd}
\end{equation}
\end{definition}
Under this formalism, we can define an unbiased monoidal category of diagrams with a product structure given by attachment. The readout map will get promoted to a functor, and compositions of diagrams are given by gluing of transfinite compositions.
\begin{definition}[category of diagrams]
For any $\alpha$-ary unbiased monoidal category $\mathscr{C}$, there is a \textit{category of diagrams} $\mathrm{dgm}(\mathscr{C})$ that shares the same objects as $\mathscr{C}$, and has morphisms that are equivalence classes of valid diagrams under the equivalence relation of isomorphism of diagrams. $\mathrm{dgm}(\mathscr{C})$ is given a strict unbiased monoidal structure via unbiased attachment of diagrams, written $\coprod$.
\end{definition}
Since isomorphisms of diagrams preserve the readout, we can treat it as an identity-on-objects functor
\begin{equation}
\mathcal{R}\colon\mathrm{dgm}(\mathscr{C})\to\mathscr{C}\qquad [\mathcal{D}]_{\simeq}\mapsto\mathcal{R}(\mathcal{D}).
\end{equation}
Moreover, it is easy to see that $\mathcal{R}$ is an unbiased monoidal functor, since 
\begin{equation}
\mathcal{R}(\coprod_{\gamma<\beta}\mathcal{D}_\gamma)=\mathcal{T}_\beta \big(\mathcal{R}(\mathcal{D}_\beta)_{\gamma<\beta}\big)
\end{equation}
and thus the data required therefor are just identities. This concludes our study of monoidal diagrams. 

\section{Conclusion}\label{conclusion}
\subsection{Summary}
Overall, we have provided a coherent theory in which infinitary monoidal diagrams can be evaluated, in the process fully constructing and exploring the notion of unbiased monoidal categories.

In Sec. \ref{cat-prelim}, we introduce the latter, starting by defining unbiased monoidal categories similarly to Ref. \citenum{unbiased-mc}, and then extending the definition to any ordinal arity. In addition, we discuss transfinite composition and prove that it commutes properly with our notion of infinitary product. Finally, we define the \textit{colimit expansion}, which both provides various examples of infinitary monoidal categories and gives rather unintuitive results regarding them.

Sec. \ref{monoidal-diagrams} defines the monoidal diagram, a novel algebraic structure which resembles a graph structure over a directed graph. After fixing the issue of $\gamma$-unresolved edges in these diagrams, we construct the readout functor which evaluates monoidal diagrams identifying them uniquely with morphisms in the underlying infinitary monoidal category. This theory is presented in such a way to respect the infinitary structure of the category; for this reason infinitary monoidal diagrams satisfy much of the same properties one derives from isotopic deformations. We conclude by promoting the class of diagrams to a category itself, making the readout an unbiased monoidal identity-on-objects functor.

\subsection{Future Work}
There are many directions in which to build upon our theory. Firstly, as has been noted in the relevant literature \cite{unbiased-mc}, there is no canonical \textit{one} extension from finitary to infinitary monoidal categories. Thus, diagrams presented here could potentially have different evaluations and satisfy different properties if a different extension is chosen. The author believes, however, that the use of unbiased monoidal categories is most natural due to the infinitary Lawvere Theory formulation discussed at the end of Sec. \ref{sec-def-unb}.

Another important choice made here was that between an algebraic theory of diagrams and a graphical theory (as it is, after all, a graphical calculus). As we chose an algebraic theory, the resultant diagrams were naturally resistant to isotopic deformation, and the difficulty in proving coherence was shunted to category theory. In contrast, working early on in a \textit{topological} formalism simplify the details of the proof. Thus, one way to generalize could be by making Joyal and Street's topological progressive plane diagrams infinitary. An important and interesting step would be proving the direct equivalence of the topological and categorical methods, thus demonstrating the true `naturality' of string diagrams for monoidal categories of unrestricted arity.

A study of the colimit-expansion presented here, the author believes, would be very fruitful. It is possible, and very likely, that many categories beyond just the semicocartesian closed ones discussed are colimit-expandable. Specifically, here we only considered satisfying the smoothness necessary for the expansion with right adjoint. But, of course, sequential cocontinuity is a much weaker condition than general cocontinuity---therefore, classifying this would provide much insight into the colimit expansion and infinitary monoidal categories.

Finally, and most apparently, one could generalize our notion of unbiased monoidal categories to the braided (and even symmetric) variants, and thereafter the diagrams that stem from it. This presents an additional challenge because of the lack of commutativity in ordinal addition. Considering our Lawvere theory approach, the isomorphism established by the braiding is mapped to by non-commutative ordinals, make it difficult to define. Furthermore, even $\mathbf{Set}$ with products would not form an example. Improving the formalism to support such notions \textit{would} be an important advancement, though, and will therefore be the future focus of the author.

\bibliography{sn-bibliography}

\end{document}